\documentclass[11pt]{article}
\usepackage{amsfonts}
\usepackage{latexsym}
\usepackage{amsmath}
\usepackage{amssymb}
\usepackage{color}
\usepackage{amsthm}
\usepackage{fullpage}
\usepackage{enumerate}
\newtheorem{theorem}{Theorem}[section]
\newtheorem{lemma}[theorem]{Lemma}

\newtheorem{proposition}[theorem]{Proposition}
\newtheorem{corollary}[theorem]{Corollary}
\newtheorem{question}[theorem]{Question}
\newtheorem{remark}[theorem]{Remark}

\theoremstyle{definition}
\newtheorem{definition}[theorem]{Definition}

\newtheorem{thmy}{Theorem}

\newtheorem*{note*}{Note}



\makeatletter
\newcommand{\subjclass}[2][1991]{%
  \let\@oldtitle\@title%
  \gdef\@title{\@oldtitle\footnotetext{#1 \emph{Mathematics subject classification.} #2}}%
}
\newcommand{\keywords}[1]{%
  \let\@@oldtitle\@title%
  \gdef\@title{\@@oldtitle\footnotetext{\emph{Key words and phrases.} #1.}}%
}
\makeatother

\DeclareMathOperator*{\esssup}{ess\,sup}
\DeclareMathOperator*{\essinf}{ess\,inf}

\begin{document}

\title{Star bodies with completely symmetric sections}

\author{Sergii Myroshnychenko, Dmitry Ryabogin, and Christos Saroglou}


\date{\today}
\maketitle
\abstract{\footnotesize We say that a star body $K$ is completely symmetric if it has centroid at the origin and its symmetry group $G$ forces any ellipsoid whose symmetry group contains $G$, to be a ball. In this short note, we prove that if all central sections of a star body $L$ are completely symmetric, then $L$ has to be a ball. A special case of our result states that if all sections of $L$ are origin symmetric and 1-symmetric, then $L$ has to be a Euclidean ball. This answers a question from \cite{R2}. Our result is a consequence of a general theorem that we establish, stating that if the restrictions in almost all equators of a real function $f$ defined on the sphere, are isotropic functions, then $f$ is constant a.e. In the last section of this note, applications, improvements and related open problems are discussed and two additional open questions from \cite{R} and \cite{R2} are answered.}
\section{Introduction}
\hspace*{1.5em}In what follows, we fix an orthonormal basis in the $n$-dimensional Euclidean space $\mathbb{R}^n$ and we denote by $\langle x,y\rangle$ the scalar product of two vectors $x,\ y$, with respect to this orthonormal basis. We denote by $S^{n-1}$ the unit sphere in $\mathbb{R}^n$, i.e. the set $\{x\in\mathbb{R}^n:|x|=1\}$, where $|x|=\sqrt{\langle x,x\rangle}$ is the length of $x$. We denote by $x^{\perp}$ the subspace which is orthogonal to $x$, i.e. $x^{\perp}=\{u\in\mathbb{R}^n: \langle x,u\rangle=0\}$. Every set of the form $S^{n-1}\cap x^{\perp}$ will be called an ``equator of $S^{n-1}$''. Denote, also, by $dx$ (or $dy$, $du$) the volume element on $S^{n-1}$ or on equators of $S^{n-1}$. The set of all isometries in $\mathbb{R}^n$ will be denoted by $ISO(n)$ and the set of all linear isometries (i.e. the orthogonal group) will be denoted by $O(n)$. Let $H$ be linear space. We denote by $SO(H)$, the set of all rotations in $H$. As usual, we set $SO(n):=SO(\mathbb{R}^n)$. A convex body in $\mathbb{R}^n$ will be a convex compact set with non-empty interior. The term ``star body'' will always refer to a set which is star-shaped with respect to the origin and has non-empty interior. A star body $K$ will be called origin symmetric (or symmetric) if $K=-K$ and $K$ will be called centrally symmetric if $K$ is a translation of $-K$. Finally, we say that $K$ is \emph{1-symmetric} if its symmetry group contains the symmetry group of a cube in $\mathbb{R}^n$.

Several important results concerning characterizations of Euclidean spaces by properties of their subspaces or quotients
have appeared in the past; we refer to \cite{Gr}, \cite{M} and \cite{Sch2}. The starting point of this manuscript is the study of the following related problem.
\begin{question}\cite{R2}\label{question-main}
Let $K$ be a convex body in $\mathbb{R}^n$, $n\geq 3$. If all orthogonal projections of $K$ are 1-symmetric, is it true that $K$ must be a Euclidean ball?
\end{question}
We will show that the answer to this question is affirmative.
It turns out that we are able to prove much more, as we provide with a rather general sufficient condition for the structure of symmetries of projections or sections of $K$ that forces $K$ to be a ball.

First we will need some definitions. Recall that the symmetry group of a set (function on $S^{n-1}$) is the group of isometries $Sym(K):=\{T\in ISO(n):TK=K\}$ (resp. $Sym(f):=\{T\in O(n):f(Tx)=f(x)\ , \ \textnormal{ for all } x\in S^{n-1}\}$).
\begin{definition}\label{def-main}

A subgroup $G$ of $ISO(n)$ will be called \emph{complete}, if every ellipsoid whose symmetry group contains $G$ is a ball.
\end{definition}
Famous examples of complete groups of isometries (see also Remarks \ref{rem-0} and \ref{rem-0'} below) are dihedral groups (in $\mathbb{R}^2$), symmetry groups of platonic solids (in $\mathbb{R}^3$), the symmetry group of the regular $n$-simplex and the symmetry group of the $n$-cube.
\begin{definition}\label{def-comp-sym-2}
A subset of $\mathbb{R}^n$ or a function $f:S^{n-1}\to\mathbb{R}$ is called \emph{completely symmetric} if its centroid is at the origin and its symmetry group is complete.
\end{definition}

We have the following result:
\begin{theorem}\label{cor-main}
Let $f:S^{n-1}\to\mathbb{R}$ be a continuous function whose restriction to every equator is completely symmetric. Then, $f$ is constant.
\end{theorem}
The following corollaries will follow immediately from Theorem \ref{cor-main}.
\begin{corollary}\label{cor-1}
Let $K$ be a star body with continuous radial function in $\mathbb{R}^n$, $n\geq 3$, whose central sections (i.e. intersections with subspaces) are all completely symmetric. Then, $K$ is an origin symmetric Euclidean ball.
\end{corollary}
For instance, if all central sections of a star body $K$ with continuous boundary are origin symmetric and 1-symmetric, then $K$ has to be an origin symmetric ball.
\begin{corollary}\label{cor-2}
Let $K$ be a convex body in $\mathbb{R}^n$, $n\geq 3$, whose orthogonal projections are all completely symmetric. Then, $K$ is an origin symmetric Euclidean ball.
\end{corollary}
\begin{corollary}\label{cor-3}
Let $K$ be a convex body in $\mathbb{R}^n$, $n\geq 3$, whose central sections are all 1-symmetric. Then, $K$ is a Euclidean ball.
\end{corollary}
\begin{corollary}\label{cor-4}
Let $K$ be a convex body in $\mathbb{R}^n$, $n\geq 3$, whose orthogonal projections are all 1-symmetric. Then, $K$ is a Euclidean ball.
\end{corollary}
We call a measurable function $f:S^{n-1}\to\mathbb{R}$ \emph{ isotropic} if the signed measure $fdx$ is isotropic. Recall that a signed Borrel measure $\mu$ on $S^{n-1}$ is called isotropic if its center of mass is at the origin and the map
$$S^{n-1}\ni x\mapsto\int_{S^{n-1}}\langle x,y\rangle^2 d\mu(y)\in\mathbb{R}$$
is constant.

Theorem \ref{cor-main} will follow easily from the following general result:
\begin{theorem}\label{thm-main}
Let $f:S^{n-1}\to\mathbb{R}$ be a measurable, bounded a.e. and even function, $n\geq 3$. If for almost every $u\in S^{n-1}$ the restriction $f|_{S^{n-1}\cap u^{\perp}}$ of $f$ to $S^{n-1}\cap u^{\perp}$ is isotropic (i.e. the restriction of $f$ to almost every equator is isotropic), then $f$ is almost everywhere equal to a constant.
\end{theorem}
Proofs of the aforementioned results will be given in Section 3. In Section 4, we establish some related results (namely, we prove a stability version of Theorem \ref{thm-main} and we answer two questions from \cite{R} and \cite{R2}) and we pose some open problems. Let us comment on the method of the proof of Theorem \ref{thm-main}: The novelty consists of the fact that we avoid using involved geometric and topological arguments, that often appear in results of Geometric Tomography (i.e. problems of determination of sets by data of their projections or sections). Instead, we make use of some integral geometric formulae together with an isoperimetric inequality and its equality cases.
\begin{remark}\label{rem-0}
In the plane, there is a very simple description of completeness. A subgroup $G$ of $ISO(2)$ is complete if and only if it contains an operator which can be written as the composition of a rotation $T$, different than $\pm Id$, and a translation. 
\end{remark}
To see this, assume for simplicity that $G\subseteq O(2)$ and assume that $G$ contains a rotation different than $\pm Id$. The reader can check that if $E$ is an origin symmetric ellipse with $TE=E$, then $E$ must be a disk. Conversely, assume that $G$ does not contain a rotation different than $\pm Id$. Since the composition of two reflections is a rotation, $G$ contains at most three elements: $Id$, $-Id$ and a reflection $R$. Then, for any origin symmetric ellipse $E$, with $RE=E$, $E$ is invariant under the action of $G$, which proves our claim.
\begin{corollary}\label{cor-3-D}
Let $f:S^2\to\mathbb{R}$ be a continuous fuction. Assume that for every $u\in S^2$, there exists $T_u\in SO(u^{\perp})\setminus\{\pm Id\}$, such that $f(T_ux)=f(x)$, for all $x\in u^{\perp}$. Then, $f$ is constant.
\end{corollary}
\noindent Proof. Immediate from Theorem \ref{cor-main} and the previous remark. $\square$
\section{Background}
In this section, we mention some basic facts about convex and star bodies that will be needed subsequently. We refer to \cite{G} and \cite{Sc} for more details and references.
Let $K$ be a star body in $\mathbb{R}^n$, that contains the origin in its interior. Its radial function $\rho_K:S^{n-1}\to \mathbb{R}_+$ of $K$ is defined by
$$\rho_K(u)=\sup\{tu:tu\in K\}\ .$$
It is true that $\rho_{K\cap u^{\perp}}=\rho_K|_{S^{n-1}\cap u^{\perp}}$,
for all $u\in S^{n-1}$. The star body $K$ is called isotropic if its centroid is at the origin and
\begin{equation}\label{eq-isotropicity-1st}
 \int_K\langle x,y\rangle^2dy=|K|^{(n+2)/n}L_K^2\ ,\qquad\textnormal{for all }x\in S^{n-1}\ ,
\end{equation}
where $L_K$ is a constant that depends only on $K$ and it is called ``the isotropic constant of K'' (see \cite{MP}). If $K$ is origin symmetric (or more generally if its centroid is at the origin), then it follows directly by integration in polar coordinates that $K$ is isotropic if and only if the function $\rho^{n+2}_K$ is isotropic. Moreover, one can easily check (see again \cite{MP}) that an equivalent statement to (\ref{eq-isotropicity-1st}) is the following:
\begin{equation}\label{eq-isotropicity-2nd}
 \int_Kx_ix_jdx=\delta_{ij}|K|^{(n+2)/n}L_K^2\ ,\qquad i,j\in\{1,\dots,n\}\ .
\end{equation}

Let $L$ be convex a convex body in $\mathbb{R}^n$. Its support function $h_L:\mathbb{R}^n\to \mathbb{R}$ is given by
$$h_L(x)=\max\{\langle x,y\rangle:y\in L\}\ .$$
Note that for $u\in S^{n-1}$, we have $h_{L|u^{\perp}}|_{S^{n-1}}=h_L|_{S^{n-1}\cap u^{\perp}}$, where $L|u^{\perp}$ denotes the orthogonal projection of $L$ onto the subspace $u^{\perp}$.
The function $h_L$ is clearly convex and positively homogeneous. Conversely, every convex and positively homogeneous function is the support function of a unique convex body.

The surface area measure of $L$ (viewed as a measure on $S^{n-1}$) is defined by $$S_L(\Omega)=\mathcal{H}^{n-1}\Big(\big\{x\in \textnormal{bd}L:\exists u\in \Omega, \textnormal{ such that }\langle x,u\rangle=h_K(u)\big\}\Big)\ ,\  \  \Omega\textnormal{ is a Borel subset of }S^{n-1}\ ,$$
where $\textnormal{bd}L$ is the boundary of $L$, $\mathcal{H}^{n-1}$ is the $(n-1)$-dimensional Haussdorff measure.
It holds that
\begin{equation}\label{eq-total-curvature}
\int_{S^{n-1}}dS_L(x)=\partial L\ ,
\end{equation}
where $\partial L$ denotes the surface area of $L$.
If $S_L$ is absolutely continuous with respect to the Lebesgue measure, its density is denoted by $f_L$ and it is usually called ``the curvature function of $L$''. 

Assume, now that the convex body $L$ contains the origin in its interior. Its polar body $L^\circ$ is defined by:
$$L^{\circ}=\{x\in\mathbb{R}^n:\langle x,y\rangle\leq 1, \textnormal{ for all }y\in L\}\ .$$
Then, $L^{\circ}$ is also a convex body that contains the origin in its interior. Moreover, $(L^{\circ})^{\circ}=L$ and
\begin{equation}\label{eq-polarity}
(L\cap u^{\perp})^{\circ}=L^{\circ}|u^{\perp}\qquad\textnormal{and}\qquad (L|u^{\perp})^{\circ}=L^{\circ}\cap u^{\perp}\ ,
\end{equation}
for all $u\in S^{n-1}$.
\section{Proofs}
\hspace*{1.5em}We start this section with the proof of our main result, Theorem \ref{thm-main}. All constants that appear will be positive constants that depend only on the dimension $n$.

Let $K$ be a symmetric star body in $\mathbb{R}^n$. The centroid body $\Gamma K$ of $K$ is the symmetric convex body whose support function is defined as follows:
$$h_{\Gamma K}(u)=\int_K|\langle x,u\rangle|dx\ .$$
Define also the quantity
\begin{equation}\label{eq-B(K)}
B(K):=\int_K\dots\int_K\textnormal{det}(x_1,\dots,x_n)^2dx_1\dots dx_n\ .
\end{equation}
We will need the following fact:
\begin{lemma}\label{lemma-centroid-density}
The surface area measure of $\Gamma K$ is absolutely continuous and its density $f_{\Gamma K}$ is given by:
$$f_{\Gamma K}(u)=c_1B(K\cap u^{\perp})\ ,\qquad u\in S^{n-1}\ .$$
\end{lemma}
\noindent Proof. 
Note that by integrating in polar coordinates, the support function of the centroid body of $K$ can be written as:
\begin{equation}\label{eq-centroid-projection-body}
h_{\Gamma K}(u)=\frac{1}{n+1}\int_{S^{n-1}}|\langle x,u\rangle|\rho_K^{n+1}(x)dx\ ,\qquad u\in S^{n-1}\ .
\end{equation}
Lemma \ref{lemma-centroid-density} is actually a reformulation of a result of Weil about the curvature function of the projection body of a convex body with absolutely continuous surface area measure. Recall the definition of the projection body $\Pi L$ of a convex body $L$:$$h_{\Pi L}(u)=\frac{1}{2}\int_{S^{n-1}}|\langle x,u\rangle| dS_L(x)\ ,\qquad u\in S^{n-1}\ . $$As Weil showed (see \cite{W}), if $f=f_L=dS_L/dx$, then for all $u\in S^{n-1}$, it holds:
\begin{equation}\label{eq-Weil}
f_{\Pi L}(u)=c_1'\int_{S^{n-1}\cap u^{\perp}}\dots\int_{S^{n-1}\cap u^{\perp}}\textnormal{det}(x_1,\dots,x_{n-1})^2f(x_1)\dots f(x_{n-1})dx_1\dots dx_{n-1}
\end{equation}
By the Minkowski Existence and Uniqueness Theorem (see again \cite{Sc}), there exists a unique (up to translation) convex body $C(K)$ with surface area measure given by
$$dS_{C(K)}=\frac{2}{n+1}\rho_K^{n+1}dx\ .$$
It is then clear that $\Gamma K=\Pi C(K)$, so by (\ref{eq-centroid-projection-body}) and (\ref{eq-Weil}), for $u\in S^{n-1}$, we get:
\begin{eqnarray*}
f_{\Gamma K}(u)&=&f_{\Pi C(K)}(u)\\
&=&c_1''\int_{S^{n-1}\cap u^{\perp}}\dots\int_{S^{n-1}\cap u^{\perp}}\textnormal{det}(x_1,\dots,x_{n-1})^2\rho_K^{n+1}(x_1)\dots \rho_K^{n+1}(x_{n-1})dx_1\dots dx_{n-1}\\
&=&c_1B(K\cap u^{\perp})\ ,
\end{eqnarray*}
where the last equality follows again by integration in polar coordinates.
$\square$\\
\\

The following easy fact (see e.g. \cite{MP}) will be also needed.
\begin{lemma}\label{lemma-legendre-isotropic}
If $K$ is isotropic, then $$\int_K|x|^2dx=c_2B(K)^{1/n}\ .$$
\end{lemma}
\noindent Proof.
It follows immediately by (\ref{eq-isotropicity-1st}) and (\ref{eq-isotropicity-2nd}) that
$$\int_K |x|^2dx=n|K|^{(n+2)/n}L_K^2\ .$$
On the other hand, expanding the determinant in the definition of $B(K)$, one easily sees:
$$B(K)=c_2'(L_K^2)^n
\ ,$$
as required. $\square$\\
\\
Proof of Theorem \ref{thm-main}: Let $c\in\mathbb{R}$. It is clear by our assumption and by the obvious fact that $c$ is isotropic in every equator of $S^{n-1}$ that $f+c$ is bounded a.e., even and isotropic in almost every equator of $S^{n-1}$. Since $f$ is bounded, we may choose $c$ to be so large that $f+c$ is positive a.e. Therefore, by replacing $f$ by $f+c$ if necessary, we may assume that $f$ is positive a.e.. Consequently, we may view the function  $f^{1/(n+1)}$ as the radial function $\rho_K$ of an origin-symmetric star body $K$. Then, by our assumption $\rho_K^{n+1}$ is isotropic in almost every equator, thus $K\cap u^{\perp}$ is isotropic, for almost every $u\in S^{n-1}$. Integrating in polar coordinates twice, one has:
$$\int_K|x|dx=\frac{1}{n+1}\int_{S^{n-1}}\rho_K^{n+1}dx=c_3\int_{S^{n-1}}\int_{S^{n-1}\cap u^{\perp}}\rho_K^{n+1}dxdu=c_4\int_{S^{n-1}}\int_{K\cap u^{\perp}}|x|^2dxdu\ .$$
Using the assumption that $K\cap u^{\perp}$ is isotropic, together with Lemma \ref{lemma-legendre-isotropic}, we get:
$$\int_K|x|dx=c_5\int_{S^{n-1}}B(K\cap u^{\perp})^{1/(n-1)}du\ ,$$
so using H\"older's inequality we obtain:
\begin{equation}\label{eq-Holder}
\int_K|x|dx\leq c_6\Big(\int_{S^{n-1}}B(K\cap u^{\perp})du\Big)^{1/(n-1)}\ .
\end{equation}
Next, using Lemma \ref{lemma-centroid-density} and (\ref{eq-total-curvature}), (\ref{eq-Holder}) immediately becomes:
\begin{equation}\label{eq-after-Holder}
\int_K|x|dx\leq c_7\Big(\int_{S^{n-1}}f_{\Gamma K}du\Big)^{1/(n-1)}=c_8\partial \big(\Gamma K\big)^{1/(n-1)}\ ,
\end{equation}
where $\partial L$ denotes the surface area of a convex body $L$. Recall the particular case of the Aleksandrov-Fenchel inequality (see e.g. \cite{Sc}):
\begin{equation}\label{eq-A-F}
\big(\partial L\big)^{1/(n-1)}\leq c_0W(L):=\int_{S^{n-1}}h_Ldx\ ,
\end{equation}with equality if and only if $L$ is a ball (note that the quantity $W(L)$ is proportional to the mean width of $L$). Applying this to (\ref{eq-after-Holder}), we deduce:
\begin{equation}\label{eq-last}
\int_K|x|dx\leq c_9W(\Gamma K)
=c_{10}\int_K|x|dx\ .
\end{equation}
Note that equality holds in both (\ref{eq-Holder}) and (\ref{eq-A-F}) (and therefore in (\ref{eq-last})) if $K$ is a ball. Thus, $c_{10}=1$. On the other hand, if $K$ is not a ball, it is well known (see again \cite{Sc} or \cite{G}) that (since $K$ is centrally symmetric), $\Gamma K$ is not a ball, so the inequality (\ref{eq-A-F}) is strict. This is a contradiction and our result follows. $\square$\\
\\
\begin{remark}\label{rem-l1}
It follows from the proof of Theorem \ref{thm-main} that if we restrict ourselves to non-negative functions (or more generally bounded from below a.e.), the assertion of the Theorem remains true if we assume $f$ to be just integrable instead of a.e. bounded.
\end{remark}
\begin{remark}
The evenness assumption in Theorem \ref{thm-main} cannot be dropped. For instance, it is well known that in any dimension $n\geq 2$, there exist convex bodies of constant width that are not Euclidean balls, i.e. convex bodies $K$ such that $h_K+h_{-K}$ is constant but $K$ is not a ball. Clearly, the support function $h_K$ of a body of constant width $K$ is isotropic in every equator, but $h_K$ is not constant if $K$ is not a ball.
\end{remark}

Before we proceed to the proof of Theorem \ref{cor-main}, we will need the following lemma:
\begin{lemma}\label{lemma-comp-sym-implies-isotropy}
Let $f:S^{n-1}\to \mathbb{R}_+$ be an integrable function. If $f$ is completely symmetric, then $f$ is isotropic.
\end{lemma}
\noindent Proof. Let $G$ be the symmetry group of $f$. Since the centroid operator intertwines isometries, it follows immediately that $G\subseteq O(n)$. Furthermore, notice that $$\int_{S^{n-1}}f(x)\langle x,u\rangle^2dx\leq \int_{S^{n-1}}f(x)|x|^2dx=\int_{S^{n-1}}f(x)dx<\infty\ ,$$for all $u\in S^{n-1}$.
Define the ellipsoid $E$, whose support function is given by:
$$h_E(u)^2=\int_{S^{n-1}}f(x)\langle x,u\rangle^2 dx\ ,\qquad u\in S^{n-1}\ .$$
Then, for any $T\in G$ and for any $u\in S^{n-1}$, we have:
\begin{eqnarray*}
h_E(u)^2&=&\int_{S^{n-1}}f(x)\langle x,u\rangle^2 dx\\
&=& \int_{S^{n-1}}f(Tx)\langle x,u\rangle^2 dx\\
&=&\int_{S^{n-1}}f(x)\langle T^{-1}x,u\rangle^2 dx\\
&=& \int_{S^{n-1}}f(x)\langle x,Tu\rangle^2 dx=h_E(Tu)^2\ .
\end{eqnarray*}
This shows that the symmetry group of $h_E$ and hence the symmetry group of $E$ contains $G$. Since $G$ is complete, it follows that $E$ is a Euclidean ball, thus $h_E$ is constant on $S^{n-1}$. This is equivalent to the fact that $f$ is isotropic. $\square$
\begin{remark}\label{rem-0'}
In view of the proof of Lemma \ref{lemma-comp-sym-implies-isotropy}, define the $O(n)$-symmetry group of a set $K$ (resp. a function $f$ on $S^{n-1}$) to be the subgroup of $O(n)$: $Sym(K)\cap O(n)$ (resp. $Sym(f)\cap O(n)$).
The following is an equivalent definition of complete symmetry: A real function defined on the sphere or a subset of $\mathbb{R}^n$ is called completely symmetric, if its $O(n)$-symmetry group is complete.
\end{remark}
To see this, notice that since the centroid operator (for sets or functions) intertwines isometries, it follows that the centroid $x$ of all sets or functions, whose $O(n)$-symmetry group is complete, has to be the origin. Indeed, if $x$ is not the origin and $G$ is the $O(n)$-symmetry group of a set or a function on the sphere, then for all $T\in G$, we must have $Tx=x$. This shows that the line $\mathbb{R}x$ is invariant under the action of $G$. Then, for any origin-symmetric ellipsoid $E$ one of whose major axis' is $\mathbb{R}x$ and its intersection with the hyperplane $x^{\perp}$ is an $(n-1)$-dimensional ball, $Sym(E)\supseteq G$, but is not necessarily a ball. Thus, the centroid of a set or function, whose $O(n)$-symmetry group is complete, is always the origin. But then, as in the proof of Lemma \ref{lemma-comp-sym-implies-isotropy}, it follows that its symmetry group is contained in $O(n)$.\\
\\
Proof of Theorem \ref{cor-main}: As in the proof of Theorem \ref{thm-main}, we may assume that $f$ is non-negative. Then, for every $p> 0$, the function $f^p$ is completely symmetric in every equator of $S^{n-1}$, thus by Lemma \ref{lemma-comp-sym-implies-isotropy}, $f(x)^p$ is isotropic in every equator of $S^{n-1}$. But then, it is immediate to check that the function $f(-x)^p$ is also isotropic in every equator of $S^{n-1}$. It follows that the function $F_p(x):=(1/2)(f(x)^p+f(-x)^p)$ is isotropic in every equator of $S^{n-1}$. Since $F_p$ is also even for all $p\geq 0$ and since $F_p$ is continuous, it follows by Theorem \ref{thm-main} that $F_p\equiv d_p$, where $d_p\geq 0$ is a constant which depends only on $p$. In particular, $f(x)+f(-x)\equiv 2d_1$. Furthermore, we get $\max\{f(x),f(-x)\}\equiv d_{\infty}:=\lim_{p\to\infty}(d_p)^{1/p}$. Thus, for every $x\in S^{n-1}$, we have $f(x)=d_{\infty}$ or $f(x)=2d_1-d_{\infty}$. This together with the continuity of $f$ prove our claim. $\square$
\\
\\
Proof of Corollaries \ref{cor-1} and \ref{cor-2}: Let $K$ be a star body with continuous radial function whose sections are completely symmetric (resp. a convex body whose projections are completely symmetric). Then, $\rho_K$ (resp. $h_K$) is continuous and completely symmetric in every equator of $S^{n-1}$. This shows that $\rho_K$ (resp. $h_K$) is constant thus $K$ is a ball. $\square$
\\
\\
Proof of Corollaries \ref{cor-3} and \ref{cor-4}: It is clear that an 1-symmetric convex body is centrally symmetric and if its center of symmetry is the origin, then it is also completely symmetric. Let $K$ be a convex body with 1-symmetric sections (resp. projections).
It is known (see \cite[Corollaries 7.1.3 and 3.1.5]{G}) that if a convex body $K$ has centrally symmetric central sections (resp. projections), then it is itself centrally symmetric. Therefore, there is a vector $x\in \mathbb{R}^n$, such that all sections (resp. projections) of $K+x$ are completely symmetric, which by Corollary \ref{cor-1} (resp. \ref{cor-2}) shows that $K+x$ is an origin symmetric ball, hence $K$ is a ball, as claimed. $\square$
\section{Further remarks and open problems}
$\S$ 1. As we have seen, the method for proving Theorem \ref{thm-main} (and therefore all its applications described here) is remarkably quick. Another advantage is that a modification of this method implies an improvement of Theorem \ref{thm-main} for positive functions, namely a stability result. To demonstrate this, let us recall that $B(K)$ (defined by (\ref{eq-B(K)})) is invariant under volume-preserving linear transformations. This, easily implies that if $K$ is a compact set with centroid at the origin, then $$\int_K|x|^2dx\geq c_2B(K)^{1/n}\ ,$$with equality if and only if $K$ is isotropic ($c_2$ is the constant from Lemma (\ref{lemma-legendre-isotropic}). It is, therefore, reasonable to define $K$ to be ``$\varepsilon$-isotropic'' (for some $\varepsilon\geq 0$), if $K$ has centroid at the origin and $\int_K|x|^2dx\leq (1+\varepsilon)c_2B(K)^{1/n}$. Similarly, for a non-negative function $f:S^{n-1}\to\mathbb{R}$, we say that $f$ is $\varepsilon$-isotropic, if
\begin{equation}\label{eq-1+epsilon-isotropic}
\int_{S^{n-1}}fdx\leq (1+\varepsilon)\overline{c}\bigg(\int_{S^{n-1}}\dots\int_{S^{n-1}}\textnormal{det}(x_1,\dots,x_n)^2f(x_1)\dots f(x_n)dx_1\dots dx_n\bigg)^{1/n} \ ,
\end{equation}
where the constant $\overline{c}$ is chosen so that there is equality in (\ref{eq-1+epsilon-isotropic}) when $f\equiv 1$ and $\varepsilon=0$. Note that a set or a function is isotropic if and only if it is 0-isotropic. Integration in polar coordinates implies immediately that a symmetric star body (or more generally a star body with centroid at the origin) $K$ is $\varepsilon$-isotropic if and only if the function $\rho_K^{n+2}$ is $\varepsilon$-isotropic. Now, let $f=\rho_K^{n+1}:S^{n-1}\to (0,\infty)$ ($K$ is as always a symmetric star body) be such that $f|_{S^{n-1}\cap u^{\perp}}$ is $\varepsilon$-isotropic, for all $u\in U$, where $U$ is a measurable subset of $S^{n-1}$, whose complement has measure less than $\delta>0$. Following the steps of the proof of Theorem \ref{thm-main}, we arrive at:
\begin{eqnarray*}
c_0W(\Gamma K)\leq(1+\varepsilon)\Big(\partial(\Gamma K)\Big)^{1/(n-1)}+\overline{c}'\delta\esssup\rho_K^{n+1}  \leq \big(\partial(\Gamma K)\big)^{1/(n-1)}+\overline{c}''(\varepsilon+\delta)\esssup\rho_K^{n+1}\ ,
\end{eqnarray*}
where $\overline{c}'$ and $\overline{c}''$ are positive constants that depend only on $n$. On the other hand, since trivially $\Gamma K$ contains the centroid body of the ball of radius $\essinf \rho_K$ and is contained in the centroid body of the ball of radius $\esssup\rho_K$ (recall that the centroid body of an origin symmetric ball is a ball), it follows from \cite[(7.124)]{Sc} that $$c_0W(\Gamma K)-(\partial (\Gamma K))^{1/(n-1)}\geq Cd_H(\Gamma K, \overline{B})\ ,$$
where $d_H$ is the Haussdorff metric, $\overline{B}$ is the ball centered at the origin with the same mean width as $\Gamma K$ and $C$ is a constant that depends only on $n$, $\esssup \rho_K^{n+1}$ and $\essinf\rho_K^{n+1}$. Since $(n+1)h_{\Gamma K}(u)=\int_{S^{n-1}}f(x)|\langle x,u\rangle|dx$ and $W(\overline{B})=W(\Gamma K)$ is proportional to $\int_{S^{n-1}}fdx$, we conclude the following:
\begin{theorem}\label{thm-stability}
There exists a function $C:(0,\infty)^2\to (0,\infty)$, with the following property: Let $\varepsilon,\delta>0$ and let $f:S^{n-1}\to\mathbb{R}$ be an even, measurable, strictly positive a.e. and bounded a.e. function, such that $f|_{S^{n-1}\cap u^{\perp}}$ is $\varepsilon$-isotropic, for all $u\in U$, for some measurable set $U\subseteq S^{n-1}$ with $|S^{n-1}\setminus U|<\delta$. Then,
$$\sup_{u\in S^{n-1}}\bigg|\int_{S^{n-1}}f(x)|\langle x,u\rangle |dx-k\int_{S^{n-1}}f(x)dx\bigg|\leq C(\essinf f,\esssup f)\cdot (\varepsilon+\delta)\ ,$$
where $k:=\int_{S^{n-1}}|x_1 |dx$.
\end{theorem}
To see that Theorem \ref{thm-stability} is actually stronger than Theorem \ref{thm-main}, recall the fact (see again \cite{Sc} or \cite{G}) that the cosine transform: $f\mapsto \int_{S^{n-1}}f(x)|\langle x,u\rangle |dx$ is injective in the class of even functions.
\\
\\
$\S$ 2. Another possible strengthening of Theorem \ref{thm-main} would be a local version of it. More precisely, it would be extremely interesting if the answer to the following question was affirmative:
\begin{question}\label{question-open-1}
Assume that for a measurable subset $U$ of $S^{n-1}$ and for an even bounded measurable function $f:S^{n-1}\to\mathbb{R}$, $f|_{S^{n-1}\cap u^{\perp}}$ is isotropic, for all $u\in U$. Is it true that $f$ is a.e. equal to a constant on the set $\bigcup_{u\in U} (S^{n-1}\cap u^{\perp})$?
\end{question}
Question \ref{question-open-1} is closely related to \cite[Problem 10]{R2}. We believe that such a result, if true, would have various applications in Geometric Tomography. To illustrate the usefulness of our conjectured local version of Theorem \ref{thm-main}, let us mention that together with Lemma \ref{lemma-old} below and an inductive argument, this would provide a quick alternative proof of a result due to the second named author \cite{R}: Let $f,\ g$ be two continuous functions on $S^{n-1}$, $n\geq 3$. If for each 2-dimensional subspace $H$, there exists $T_H\in SO(H)$, such that $f(x)=g(T_Hx)$, for all $x\in S^{n-1}\cap H$, then $f(x)=g(x)$, for all $x\in S^{n-1}$ or $f(x)=g(-x)$, for all $x\in S^{n-1}$.
\begin{lemma} \label{lemma-old}\cite{R}
Let $f,\ g$ be two continuous functions on $S^2$, such that for every $u\in S^2$, $f(x)=g(x)$, for all $x\in u^{\perp}$ or $f(x)=g(-x)$, for all $x\in u^{\perp}$. Then, $f(x)=g(x)$, for all $x\in S^{2}$ or $f(x)=g(-x)$, for all $x\in S^{2}$.
\end{lemma}
\begin{remark}\label{rem-1}
Our methods allow us to drop the continuity assumption in Lemma \ref{lemma-old} (a fact which was conjectured in \cite{R}). More precisely, we can show the following: If $f$, $g$ are $L^1$ functions on $S^{n-1}$, such that for almost every $u\in S^{n-1}$, $f(x)=g(x)$, for almost every $x\in S^{n-1}\cap u^{\perp}$ or $f(x)=g(-x)$, for almost every $x\in S^{n-1}\cap u^{\perp}$, then $f(x)=g(x)$, a.e. in $S^{n-1}$ or $f(x)=g(-x)$, a.e. in $S^{n-1}$. The proof is an immediate consequence of a well know integral-geometric formula due to Bussemann (see also \cite{S-W} for generalizations and references): For any non-negative functions $F_1,\dots,F_{n-1}$ on the sphere $S^{n-1}$, we have:
\begin{eqnarray*}&&\int_{S^{n-1}}F_1dx\dots\int_{S^{n-1}}F_{n-1}dx\\
&=&c\int_{S^{n-1}}\int_{S^{n-1}\cap u^{\perp}}\dots\int_{S^{n-1}\cap u^{\perp}}|det(x_1,\dots,x_{n-1})|F(x_1)\dots F(x_n)dx_1\dots dx_{n-1}du\ ,
\end{eqnarray*}
as long as both parts exist. Here as always, $c>0$ is a constant that depends only on the dimension. Indeed, setting $F_1(x):=(f(x)-g(x))^2$, $F_2(x)=\dots =F_{n-1}(x)=(f(x)-g(-x))^2$, our assumption gives $$\int_{S^{n-1}}(f(x)-g(x))^2dx\bigg(\int_{S^{n-1}}(f(x)-g(-x))^2dx\bigg)^{n-2}=0\ ,$$which shows that $f(x)=g(x)$ a.e or $f(x)=g(-x)$ a.e. as claimed.
\end{remark}
\textit{}\\
$\S$ 3. Recall that a subset of $\mathbb{R}^n$ is called \emph{unconditional} if it is symmetric with respect to the hyperplanes $e_1^{\perp},\dots,e_{n}^{\perp}$, for some orthonormal basis $\{e_1,\dots,e_n\}$ of $\mathbb{R}^n$.
In view of Question \ref{question-main}, observe that ellipsoids and bodies of revolution have the following property: Every section (projection) is unconditional. One might ask if these are the only examples of convex bodies with this property. The answer to this question is negative as the following shows.

\begin{proposition}\label{prop-counterexample}
There exists an origin symmetric convex body in $\mathbb{R}^n$, $n\geq 3$, which is not an ellipsoid or a body of revolution and whose sections (projections) are all unconditional.
\end{proposition}
\noindent Proof. As (\ref{eq-polarity}) shows, it suffices to find a counterexample for projections, since its polar would give a counterexample for sections. It is well known that the Minkowski sum of two ellipsoids is not always an ellipsoid. It follows that the sum of an ellipsoid and a ball is not always an ellipsoid (otherwise, we would be able to apply linear transformations to show that the sum of any two ellipsoids is an ellipsoid). It follows by continuity (since the family of ellipsoids which are not bodies of revolution is dense to the family of all ellipsoids, in the sense of the Haussdorff metric) that we can find an ellipsoid $E$ which is not a body of revolution and a ball $B$ such that $E+B$ is not an ellipsoid. On the other hand, $E+B$ is clearly not a body of revolution but all its projections are unconditional. $\square$
\\

It should be remarked that Proposition \ref{prop-counterexample} gives a negative answer to Problem 13 from \cite{R2}.
Nevertheless, the following questions remain open:
\begin{question}\label{question-open-2}
Is it true that if all projections of a convex body $K$ in $\mathbb{R}^n$, $n\geq 3$, are translations of completely symmetric sets, then $K$ is a ball?
\end{question}
\begin{remark}\label{rem-2}
It should be noted that if all projections of $K$ are translations of completely symmetric sets, then the restriction of the (even) function $h_K(x)+h_K(-x)=h_K(x)+h_{-K}(x)$ on every equator of $S^{n-1}$ is completely symmetric. Thus, Theorem \ref{thm-main} immediately implies that $h_K+h_{-K}$ is constant everywhere on $S^{n-1}$, i.e. $K$ is a convex body of constant width.
\end{remark}
\begin{question}\label{question-open-3}
Is it true that if all projections in $\mathbb{R}^n$, $n\geq 3$, of a convex body $K$ are linear (or more generally affine) images of completely symmetric sets (e.g. all projections are linear images of 1-symmetric sets), then $K$ is an origin-symmetric ellipsoid?
\end{question}
Variants about sections of star bodies (or equivalently restrictions of functions) are also interesting. In this connection, see again \cite{Gr} or \cite{M}.

\vspace{0.8 cm}

\noindent Sergii Myroshnychenko\\
Department of Mathematical Sciences\\
Kent State University\\
Kent, OH 44242, USA \\
E-mail address: smyroshn@kent.edu

\vspace{0.5 cm}

\noindent Dmitry Ryabogin \\
Department of Mathematical Sciences\\
Kent State University\\
Kent, OH 44242, USA \\
E-mail address: ryabogin@math.kent.edu

\vspace{0.5 cm}

\noindent Christos Saroglou \\
Department of Mathematical Sciences\\
Kent State University\\
Kent, OH 44242, USA \\
E-mail address: csaroglo@kent.edu \ \&\ christos.saroglou@gmail.com

\end{document}